\input amstex
\pageno=1

\magnification=1200
\loadmsbm
\loadmsam
\loadeufm
\input amssym
\UseAMSsymbols

\TagsOnRight

\hsize172 true mm
\vsize216 true mm
\voffset=20 true mm
\hoffset=0 true mm
\baselineskip 9.5 true mm plus0.4 true mm minus0.4 true  mm

\def\P{\partial}

\def\X{\xi}

\vskip 0.4cm
\noindent

\centerline{Integrable systems and effectivisation }
\centerline{of Riemann theorem about domains of the complex plane}

\vskip 1cm
\centerline{S.M.Natanzon}

1. Consider a closed analytic curve $\gamma$ in the complex plane
and denote by $D_+$ and $D_-$ the interior and exterior domains
with respect to the curve. The point $z=0$ is assumed to be in $D_+$.
Then according to Riemann theorem there exists a
function $w(z)=\frac 1r z+\sum\limits_{j=0}^\infty p_j z^{-j}$,
mapping $D_-$ to the exterior of the unit disk
$\{w\in\Bbb C\vert\vert w\vert>1\}$. It is follow from [1] that
this function is described by formula
$\log w=\log z-\P_{t_0}
(\frac 12\P_{t_0}+\sum\limits_{k\geqslant 1}\frac{z^{-k}}{k}
\P_{t_k})v$, where $v=v(t_0, t_1, \bar t_1, t_2, \bar t_2,...)$
is a function from infinite number of variables
$$t_0=\frac 1\pi\int\limits_{D_+}d^2z, \quad
t_k=\frac {1}{\pi k} \int\limits_{D_-}z^{-k}d^2 z.$$
Moreover, this function satisfies the dispersionless Hirota
equation for $2D$ Toda lattice hierarchy.
$$(z-\X)e^{D(z)D(\X)v}=z e^{-\P_{t_0}D(z)v}-\X e^{-\P_{t_0}D(\X)v}.
\tag 1$$
$$(\bar z-\bar\X)e^{\bar D(\bar z)\bar D(\bar\X)v}=
\bar z e^{-\P_{t_0}\bar D(\bar z)v}-\bar\X e^{-\P_{t_0}\bar D(\bar\X)v}.
\tag 2$$
$$1- e^{-D(z)\bar D(\bar\X)v}=\frac{1}{z\bar\X} e^{\P_{t_0}
(\P_{t_0}+D(z)+\bar D(\bar\X))v}, \tag 3$$
where
$$D(z)=\sum_{k\geqslant 1}\frac{z^{-k}}{k}\P_{t_k}, \quad
\bar D(\bar z)=\sum_{k\geqslant 1}
\frac{\bar z^{-k}}{k}\P_{\bar t_k}.$$

Thus for an effectivisation of Riemann theorem it is
sufficiently to find a representation of $v$ in the form 
of Taylor series 
$$v=\sum N\bigl(i_0\bigl\vert i_1,...,i_k\bigl\vert
\bar i_1,...,\bar i_{\bar k}\bigr)t_0
t_{i_1},...,t_k \bar t_{\bar i_1},...,\bar t_{\bar i_{\bar k}}.$$
The numbers $N(i_0\bigl\vert i_1,...,i_k\bigl\vert\bar i_1,...,
\bar i_{\bar k})$ for $i_\alpha, \bar i_\beta\leqslant 2$ is 
found in [ 1 ]. In this paper we find some recurrence relations, 
that give a possible to find all
$N(i_0\bigl\vert i_1,...,i_k\vert\bar i_1,...,\bar i_{\bar k})$.
For this we find some formulas for reconstruction any solutions of 
(1) -- (3) via arbitrary Cauchy data
$\frac{\P^2 v}{\P t_0\P t_i}\Bigl\vert_{t_0}, 
\frac{\P^2 v}{\P t_0\P\bar t_i}\Bigl\vert_{t_0}$, 
$\frac{\P^2 v}{dt^2}\Bigl\vert_{t_0}$ (here $F\Bigl\vert_{t_0}=
F\Bigl\vert_{t_1=\bar t_1=t_2=\bar t_2=\dotsb=0}(t_0)$) and
find the values of these data for our $v$.

2. Let us put $\P_i=\P_{t_i}=\frac{\P}{\P t_i}$ and
 $\bar\P_i=\P_{\bar t_i}=\frac{\P}{\P \bar t_i}$.

{\bf Lemma 1.} {\it
$z-\sum\limits^\infty_{j=1}\frac {1}{j}z^{-j}\P_1\P_jv=
ze^{-\P_0D(z)}v$ \  and $$\P_1\P_j v=\sum\limits^\infty_{m=1}
\frac {(-1)^{m+1}}{m!}
\sum\limits_{k_1+\dotsb+k_m=j+1}
\frac{j}{k_1\dotsb k_m}\P_0\P_{k_1} v
\dotsb\P_0\P_{k_m}v.$$}

Proof: According (1)  $(z-\X) e^{D(z)D(\X)v}=
(z-\X)(1+(D(z) D(\X)v)+\frac 12(D(z) D(\X)v)^2+\dotsb)=
(z-\X)(1+z^{-1}\X^{-1}\P^2_1 v+z^{-1}\sum\limits_{j=2}^\infty
\frac 1j \X^{-j}\P_1\P_jv+\X^{-1}\sum\limits_{j=2}^\infty
\frac 1j z^{-j}\P_1\P_j v+z^{-2}\X^{-2}F)=(z-\X)+\X^{-1}\P^2_1 v-
z^{-1}\P^2_1 v+\sum\limits_{j=2}^\infty\frac 1j\X^{-j}\P_1\P_j v-
\sum\limits_{j=2}^\infty\frac 1j z^{-j}\P_1\P_j v+
z^{-1}\X^{-1}F.$

On the other hand according (1) the function
$(z-\X)e^{D(z) D(\X)v}$ is a sum of two functions
$f_1(z)+f_2(\X)$. Thus 
$F=0$ and $ze^{-\P_0 D(z)}v= z-\sum\limits_{j=1}^\infty\frac 1j
z^{-j}\P_1\P_j v$. Therefore,
$$\sum\limits_{j=1}^\infty
\frac 1j z^{-(j+1)}\P_1\P_j v=1-e^{-\P_0 D(z)}v=
1-\bigl(1+\sum\limits_{m=1}^\infty\frac {(-\P_0 D(z)v)^m}{m!}\bigr)=$$
$$-\sum\limits_{m=1}^\infty\frac{(-1)^m}{m!}
\Bigl(\sum\limits_{k=1}^\infty\frac{z^{-k}}{k}\P_0\P_k v\Bigr)^m=$$
$$=-\sum\limits_{m=1}^\infty\frac{(-1)^m}{m!}
\Bigl(\sum\limits_{n=1}^\infty z^{-n}
\sum\limits_{k_1+\dotsb +k_m=n}\frac{1}{k_1\dotsc k_m}
\P_0\P_{k_1} v\dotsb \P_0 \P_{k_m}v \Bigr)=$$
$$=-\sum\limits_{n=1}^\infty z^{-n}\Bigl(\sum\limits_{m=1}^\infty
\frac{(-1)^m}{m!}\sum\limits_{k_1+\dotsb +k_m=n}
\frac{1}{k_1\dotsc k_m}\P_0\P_{k_1}v\dotsb \P_0\P_{k_m} v\Bigr).$$
Thus
$$\frac 1j \P_1 \P_j v=-\sum\limits_{m=1}^\infty \frac{(-1)^m}
{m!}\sum\limits_{k_1+\dotsb +k_m=j+1}\frac{1}{k_1\dotsc k_m}
\P_0\P_{k_1}v\dotsb \P_0\P_{k_m} v. \square$$

{\bf Lemma 2.} {\it $\P_i\P_j v= \sum
\limits_{m=1}^\infty
\sum\limits_{\matrix s_1+\dotsb +s_m=j+1 \\ s_i>1\endmatrix}
\frac{(-1)^{m+1}}{m}
\frac{ij}{(s_1-1)\dotsb(s_m-1)} P_{ij}(s_1-1,...,s_m-1)\cdot$
\linebreak $\cdot\P_1 \P_{s_1-1} v\dotsb$ $\P_1\P_{s_m-1} v$,
where $P_{ij}(s_1-1,...,s_m-1)$ is the number 
of representations of
$\{i=i_1+\dotsb+i_m\vert 1\leqslant i_k\leqslant s_k-1,\ k=1,
\dotsb, m\}$ of number $i$.}

Proof: According to lemma 1 and equation (1)
$(z-\X)e^{D(z) D(\X) v}=$ $$=z-\sum\limits_{j=1}^\infty\frac 1j
z^{-j}\P_1\P_j v-(\X-\sum\limits_{j=1}^\infty\frac 1j
\X^{-j}\P_1\P_j v)=(z-\X)-\sum\limits_{j=1}^\infty\frac 1j
(z^{-j}-\X^{-j})\P_1\P_jv.$$ Thus,
$$e^{D(z) D(\X)}=1+z^{-1}\X^{-1}\sum\limits_{j=1}^\infty\frac 1j
\frac{(z^{-j}-\X^{-j})}{(z^{-1}-\X^{-1})}\P_1\P_j v=$$
$$=1+z^{-1}\X^{-1}\sum\limits_{j=1}^\infty\frac 1j
\big(\sum\limits_{\matrix s+t=j-1 \\ s, t\geqslant 0 \endmatrix}
z^{-s}\X^{-t})\P_1\P_j v=$$
$$=1+\sum\limits_{j=1}^\infty\frac 1j(\sum
\limits_{\matrix s+t=j+1 \\ s,t\geqslant 1\endmatrix}
z^{-s}\X^{-t})\P_1\P_j v.$$ Therefore,
$$D(z)D(\X)v=\sum\limits_{m=1}^\infty\frac{(-1)^{m+1}}{m}
\bigl(\sum\limits_{n=1}^\infty\bigl(\sum\limits_{\matrix s+t=n+1
\\s,t\geqslant 1 \endmatrix}
z^{-s}\X^{-t}\bigr)\frac 1n \P_1\P_n v\bigr)^m=$$
$$=\sum\limits_{j=1}^\infty\frac{(-1)^{m+1}}{m}
\sum\limits_{i,j\geqslant 1} z^{-i}\X^{-j}\cdot$$
$$\Bigl(\sum\limits_{\matrix i_1+\dotsb+i_m=i \\ j_1+\dotsb j_m=j
\\ i_k, j_k\geqslant 1\endmatrix} \frac{1}{i_1+j_1-1}
\P_1\P_{i_1+j_1-1}v\dotsb\frac{1}{i_m+j_m-1}\P_1\P_{i_m+j_m-1}v
\Bigr),$$ that is  
$$\P_i\P_j v= \sum\limits_{m=1}^\infty
\sum\limits_{s_1+\dotsb+s_m=i+j}\frac {(-1)^{m+1}}{m}
\frac{ij}{(s_1-1)\dotsb (s_m-1)}\cdot$$
$$\cdot P_{ij}(s_1-1,...,s_m-1)
\P_1\P_{s_1-1} v\dotsb\P_1\P_{s_m-1} v. \square$$

{\bf Remark.} The equations
$$\P_i\P_j v= \sum\limits_{m=1}^\infty
\sum\limits_{s_1+\dotsb+s_m=i+j}\frac
{(-1)^{m+1}}{m}\frac{ij}{(s_1-1)\dotsb (s_m-1)}\cdot$$
$$\cdot P_{ij}(s_1-1,...,s_m-1)
\P_1\P_{s_1-1} v\dotsb\P_1\P_{s_m-1} v$$
describe the dispersionless limit of KP equation.
Some other description of this hierarchy is presented
in [ 2 ]. A comparison of these descriptions gives some
nontrivial combinatorial identity $P_{ij}$.

{\bf Lemma 3.} {\it $\P_i\P_j v= \sum\limits_{m=1}^\infty
\sum\limits_{p_1+\dotsb+p_m=j+i} \frac{ij}{p_1\dotsb p_m}
T_{ij}(p_1\dotsb p_m) \P_0 \P_{p_1} v\dotsb \P_0\P_{p_m} v$,
where \linebreak $T_{ij}(p_1...p_m)=\sum\limits_
{\matrix n_1+\dotsb+n_k=m \\ n_i>0\endmatrix}
\frac{(-1)^{m+1}}{k}
\frac {1}{n_1!\dotsb n_k!} P_{ij}(p_1+\dotsb+p_{q_1}-1,...,
p_{q_{k-1}+1}+\dotsb+p_{q_k}-1)$, and $q_j=\sum\limits^j_{i=1}n_i$.}

Proof: According to lemmas 1 and 2
$$\P_i\P_j v= \sum\limits_{m=1}^\infty
\sum\limits_{s_1+\dotsb+s_m=j+i} \frac{(-1)^{m+1}}{m}
\frac{ij}{(s_1-1)\dotsb (s_m-1)}
P_{ij}(s_1-1,\dotsb, s_m-1)\cdot$$
$$\cdot\P_1 \P_{s_1-1} v\dotsb \P_1\P_{s_m-1} v)=
\sum\limits_{m=1}^\infty
\sum\limits_{s_1+\dotsb+s_m=j+i}\frac{(-1)^{m+1}}{m}
\frac{ij}{(s_1-1)\dotsb (s_m-1)}\cdot$$ $$\cdot P_{ij}(s_1-1,\dotsb,
s_m-1)\Bigl(\sum^\infty_{n_1=1}\sum\limits_{p_1
+\dotsb+p_{n_1}=s_1}
\frac{(-1)^{n_1+1}}{n_1!}\frac{s_1-1}{p_1\dotsb p_{n_1}}
\P_0 \P_{p_1} v\dotsb \P_0\P_{p_{n_1}}v\Bigl)\dotsb$$
$$\dots\Bigl(\sum^\infty_{n_m=1}
\sum\limits_{p_1+\dotsb +p_{n_m}=s_m}\frac{(-1)^{n_m+1}}{n_m!}
\frac{s_m-1}{p_1\dotsb p_{n_m}}
\P_0\P_{p_1} v\dotsb \P_0\P_{p_{n_m}} v\Bigl)=$$
$$=\sum\limits_{m=1}^\infty
\sum\limits_{p_1+\dotsb+p_m=j+i} \frac{ij}{p_1\dotsb p_m}
T_{ij}(p_1\dotsb p_m) \P_0 \P_{p_1} v\dotsb \P_0\P_{p_m}v. \square$$

Using the induction, we get from lemma 3

{\bf Lemma 4.} {\it
$$\P_{i_1}\P_{i_2}\dotsb \P_{i_k}v=\sum\limits_{m=1}^\infty
\Bigl(\sum\limits_{\matrix s_1+\dotsb+s_m=i_1+\dotsb +i_k \\ \ell_1
+\dotsb +\ell_m=m+k-2 \\ s_j, \ell_j\geqslant 1\endmatrix}
\frac{i_1\dotsb i_k}{s_1\dotsb s_m}\cdot$$
$$\cdot T_{i_1\dotsb i_k}
\pmatrix s_1\dotsb s_m \\ \ell_1\dotsb \ell_m\endpmatrix
\P^{\ell_1}_0\P_{s_1} v\dotsb \P^{\ell_m}_0\P_{s_m} v\Bigr),$$
$$ \text{where}\quad
T_{i_1, i_2}\pmatrix s_1\dotsb s_m \\ \ell_1 \dotsb \ell_m\endpmatrix=
\cases T_{i_1 i_2}(s_1\dotsb s_m),& \ \text{if}\ \ell_1=\dotsb=
\ell_m=1 \\ 0 & \ \text{in another cases}\endcases,$$
$$T_{i_1\dotsb i_k}\pmatrix s_1\dotsb s_m \\ \ell_1\dotsb \ell_m
\endpmatrix=
\sum\limits_{\matrix 1\leqslant i\leqslant j\leqslant m \\ s, \ell>0
\endmatrix} T_{i_1\dotsb i_{k-1}}
\biggl( {\matrix s_1\dotsb s_{i-1} \\ \ell_1\dotsb \ell_{i-1}\endmatrix}
\pmatrix s \\ \ell \endpmatrix {\matrix s_{j+1}\dotsb s_m\\ \ell_{j+1}
\dotsb \ell_m\endmatrix} \biggr)\cdot$$
$$\cdot T_{s,i_k}(s_i, s_{i+1},...,s_j)
\frac{\ell!}{(\ell_i-1)!\dotsb (\ell_j-1)!},$$
$$\text{and} \quad s=s_i+s_{i+1}+\dotsb +s_j-i_k, \quad
\ell=(\ell_i-1)+\dotsb+(\ell_j-1). \square $$}

3. Define now the Cauchy data for $v$.

{\bf Lemma 5.} {\it
$$\P_0 v\vert_{t_0}=-t_0+t_0\ln t_0\quad
\P_k v\vert_{t_0}=0\quad \text{for}\quad k>0.$$}

Proof: Put $\gamma=\{z\in\Bbb C\vert\vert z\vert=R\}$. 
The Schwarz function for $\gamma$ is
$S(z)=R^2z^{-1}$. Thus according to 
[ 1 (26) ] $$t_0=\frac 1\pi\int\limits_{\vert z\vert\leqslant R}
d^2z=\frac{\pi R^2}{\pi}=R^2.$$
$$t_k=\frac{1}{2\pi i k}\int\limits_{\vert z\vert=R}
z^{-k} R^2 z^{-1} dz=\frac{1}{2\pi i k} R^2 \oint
\limits_\gamma z^{-(k+1)} dz=0\quad \text{for}\quad k>0$$
Therefore, according to [1 (8), (26)]
$$\P_0 v\Bigl\vert_{t_0}=\frac 2\pi\int\limits_{\vert z\vert\leqslant R}
\ln\vert z\vert d^2 z=\frac 2\pi\int\limits_0^R dr
\int\limits_0^{2\pi} d\varphi\cdot r\cdot\ln\vert r\vert=$$
$$=\frac{2\cdot 2\pi}{2\pi}\int\limits_0^R\ln\vert r\vert
dr^2=2(r^2\ln r\Bigl|_0^R-\int\limits_0^R rdr)=$$
$$=r^2\ln r^2\Bigl\vert_0^R-
2\frac 12 r^2\Bigl\vert_0^R=R^2\ln R^2-R^2=
-t_0+t_0\ln t_0.$$
$$\P_k v\Bigl\vert_{t_0}=\frac {1}{2\pi i}\oint\limits_{\vert z\vert= R}
z^kR^2z^{-1} dz=\frac {R^2}{2\pi i}\oint\limits_\gamma
z^{k-1} dz=0\quad \text{for} \quad k>0. \square$$

{\bf Lemma 6.} {\it
$$\P_i\bar\P_j v\vert_{t_0}=\cases 0 &\ \text{for}\ i\ne j,
\\ i t^i_0 &\ \text{for}\ i=j. \endcases$$}

Proof: It is follow from lemma 5 that
$\P_0\P_k v\bigl\vert_{t_0}=0$ for $k>0$ and
$\P^2_0 v\vert_{t_0}=\ln t_0$. Moreover according to (3)
$$1-e^{-D(z)\bar D(\bar\X) v}=z^{-1}\X^{-1} e^{\P_0(\P_0+D(z)+
\bar D(\bar\X))v}.$$ Thus
$$e^{\P_0(\P_0+D(z)+\bar D(\bar\X))v}\Bigl\vert_{t_0}=t_0$$
and
$$-D(z)\bar D(\bar\X) v\bigl\vert_{t_0}=\ln(1-z^{-1}\bar\X^{-1}
t_0)=-\sum\limits_{k=1}^\infty k z^{-k}\X^{-k} t_0^k.$$
Therefore $\P_i\P_jv\bigl\vert_{t_0}=0$ for $i\ne j$ and 
$\P^2_i v\bigl\vert_{t_0}=i t^i_0$ $\square$.

{\bf Lemma 7.} {\it
$$\P_{i}\bar\P_{i_1}\dotsb\bar\P_{i_k}v\Bigl\vert_{t_0}=
\bar\P_i\P_{i_1}\dotsb \P_{i_k} v\bigl\vert_{t_0}=
\cases 0, &\ \text{if} \quad i_1
+\dotsb+i_k\ne i \\ i_1\dotsb i_k \frac{i!}{(i-k+1)!} t_0^{i-k+1}, &
\ \text{if}\ i=i_1+\dotsb +i_k \endcases.$$}

Proof: The differentials $\P$ and $\bar\P$ occur in
(1) -- (3)symmetrically. This gives the first equality.
Moreover according to lemmas 4 and 5,
$$\P_{i_1}\P_{i_2}\dotsb \P_{i_k} v=
\frac{i_1\dotsb i_k}{i}T_{i_1\dotsb i_k}{\pmatrix i
\\ k - 1\endpmatrix}\P_0^{k-1}\P_{i} v+$$
$$+\sum\limits_{m=2}^\infty
\Bigl(\sum\limits_{\matrix s_1+\dotsb+s_m=i_1+\dotsb +i_k \\ \ell_1
+\dotsb \ell_m=m+k-2 \\ s_j, t_j\geqslant 1\endmatrix}
\frac{i_1\dotsb i_k}{s_1\dotsb s_m} T_{i_1\dotsb i_k}
\pmatrix s_1\dotsb s_m \\ \ell_1\dotsb \ell_m\endpmatrix\cdot$$
$$\cdot\P^{\ell_1}_0\P_{s_1} v\dotsb \P^{\ell_m}_0\P_{s_m} v
\Bigr)=$$
$$=\frac{i_1\dotsb i_k}{i}T_{i_1\dotsb i_k}{\pmatrix i
\\ k - 1\endpmatrix} \P_0^{k-1}\P_{i} v
=\frac{i_1\dotsb i_k}{i}\P_0^{k-1}\P_i v,$$ 
where $i=i_1+\dotsb+i_k$.
This equality and lemma 6 give the second equality in the
affirmation of lemma 7. $\square$

{\bf Lemma 8.} {\it
$$\P_{i_1}\dotsb \P_{i_k}\bar\P_{\bar i_1}\dotsb
\bar\P_{\bar i_{\bar k}} v\Bigl\vert_{t_0}=
\sum_{i=1}^\infty N_i(i_1\dotsb i_k\bigl\vert
\bar i_1\dotsb \bar i_{\bar k})
t_0^{i-(k+\bar k)+2},\ 
\text{where} \
N_i\bigl(i_1\dotsb i_k\bigl\vert \bar i_1\dotsb \bar i_{\bar k}
\bigr)=0$$ $ \text{if}\ \sum\limits_{j=1}^{k} i_j\ne i$
or $\sum\limits_{j=1}^{\bar k}\bar i_j\ne i$. In opposite case,
$N_i\bigl(i_1\dotsb i_{i_k}\bigl\vert \bar i_1,\dotsb
\bar i_{\bar k}\bigr)=$ 
$$=\sum\limits_{m=1}^\infty
\Bigl(\sum\limits_{\matrix s_1+\dotsb+s_m=i_1+\dotsb +i_k \\ \ell_1
+\dotsb \ell_m=m+k-2 \\ s_j, \ell_j\geqslant 1\endmatrix}
i_1\dotsb i_k\bar i_k\dotsb \bar i_{\bar k}\cdot$$
$$\cdot T_{i_1\dotsb i_k}
\pmatrix s_1\dotsb s_m \\ \ell_1\dotsb \ell_m\endpmatrix
S_{\bar i_1\dotsb \bar i_{\bar k}}\Bigl(s_1,...,s_m\Bigr)\Bigr),$$
where  
$$S_{\bar i_1\dotsb \bar i_{\bar k}}(s_1,...,s_m)=
\sum\limits_{n_1+\dotsb+n_m=k} \sum
\frac{(s_1-1)!\dotsb (s_m-1)!}{(s_1-n_1+1-\ell_1)!\dotsb
(s_m-n_m+1-\ell_m)!}$$
and second sum is the sum by all partitions of the set
$\{\bar i_1,...,\bar i_{\bar k}\}$ on subsets
\linebreak $\{j^p_1,...,j^p_{n_p}\}$ $(p=1,...,m)$ such that
$\sum\limits^{n_p}_{\alpha=1}j_\alpha^p=s_p$. $\square$}

Proof: According to lemma 4
$$\P_{i_1}\dotsb\P_{i_k}
\bar\P_{\bar i_1}\dotsb \bar\P_{\bar i_{\bar k}}v=
\bar\P_{\bar i_{\bar i}}\dotsb \bar\P_{\bar i_{\bar k}}
\Bigl(\sum\limits_{m=1}^\infty
\Bigl(\sum\limits_{\matrix s_1+\dotsb+s_m=i_1+\dotsb +i_k \\ \ell_1
+\dotsb +\ell_m=m+k-2 \\ s_j, \ell_j\geqslant 1\endmatrix}
\frac{i_1\dotsb i_k}{s_1\dotsb s_m}\cdot$$
$$\cdot T_{i_1\dotsb i_k}
\pmatrix s_1\dotsb s_m \\ \ell_1\dotsb \ell_m\endpmatrix
\P^{\ell_1}_0\P_{s_1} v\dotsb \P^{\ell_m}_0\P_{s_m} v\Bigr)\Bigl)=$$
$$=\sum\Bigl(\sum\limits_{m=1}^\infty
\Bigl(\sum\limits_{\matrix s_1+\dotsb+s_m=i_1+\dotsb +i_k \\ \ell_1
+\dotsb +\ell_m=m+k-2 \\ s_j, \ell_j\geqslant 1\endmatrix}
\frac{i_1\dotsb i_k}{s_1\dotsb s_m}
T_{i_1\dotsb i_k}\cdot$$
$$\cdot\pmatrix s_1\dotsb s_m \\ \ell_1\dotsb \ell_m\endpmatrix
\P^{\ell_1}_0 \P_{s_1}
\bar\P_{\bar j_1}^1\dotsb \bar\P_{\bar j_{n_1}}^1
 v\dotsb \P^{\ell_m}_0 \P_{s_m}
\bar\P_{\bar j_1}^m\dotsb \bar\P_{\bar j_{n_m}}^m
 v\Bigr)\Bigl),$$
where the first sum is taken by all partitions of 
$\{\bar i_1,...,\bar i_{\bar k}\}$ on subset
$\{(j^1_i,...,j^1_{n_1})...(j^m_1,...,$ \linebreak
$j^m_{n_m})\}$.
According to lemma 7 this gives the affirmation of
lemma 8 $\square$.

Lemma 8 gives

{\bf Theorem.} {\it
$$v=\frac 12 t^2_0\log t_0-\frac 34 t_0+
\sum\limits_{\matrix i_1<\dotsb<i_k \\
\bar i_1<\dotsb \bar i_{\bar k}
\endmatrix}\sum_{n_j, \bar n_j=1}^{\infty}
\frac{1}{n_1!\dotsb n_k!\bar n_1!\dotsb\bar n_{\bar k}!}
N_i\Bigl({\matrix i_1...i_k \\ n_1...n_k\endmatrix}\Bigl\vert
\matrix \bar i_1...\bar i_{\bar k} \\ \bar n_1 ... \bar n_{\bar k}
\endmatrix\Bigr)\cdot$$ $$\cdot t_0^{i-(\sum\limits_{i=1}^k n_i+
\sum\limits_{i=1}^{\bar k}\bar n_i)+2}
t_{i_1}^{n_1}\dotsb t_{i_k}^{n_k} \bar t_{\bar i_1}^{\bar n_1}
\dotsb \bar t_{\bar i_{\bar k}}^{\bar n_{\bar k}},$$
where
$$N_i\Bigl({\matrix i_1...i_k \\ n_1...n_k\endmatrix}\Bigl\vert
\matrix \bar i_1...\bar i_{\bar k} \\ \bar n_1 ... \bar n_{\bar k}
\endmatrix\Bigr) =0,$$ if
$i\ne\sum\limits_{j=1}^k n_j i_j$ or $i\ne\sum\limits_{j=1}^{\bar k}
\bar n_j\bar i_j.$ In opposite case
$$N_i\Bigl({\matrix i_1...i_k \\ n_1...n_k\endmatrix}\Bigl\vert
{\matrix \bar i_1...\bar i_{\bar k} \\ \vert n_1 ... \bar n_{\bar k}
\endmatrix}\Bigr) =N_i(i_1...i_1i_2...i_2...i_k...i_k\bigl\vert
\bar i_1...\bar i_1\bar i_2...\bar i_2...\bar i_{\bar k}...
\bar i_{\bar k}),$$ where in the last parentheses any
$i_j$ (respectively any $\bar i_j$ ) occurs $n_j$ (respectively
$\bar n_j$) cases.}

\centerline {Reference.}

1. I.K.Kostov, I.Krichever, M.Mineev--Weinstein, P.B.Wiegmann, A.Zabrodin,
\linebreak $\tau$-function for anatytic curves, arXiv:hep-th/0005259.

2. S.M.Natanzon, Formulas for $A_n$ and $B_n$--solutions of
WDVV equations, Journal of Geometry and Physics,  (2001).

\end